\documentstyle[11pt,epic]{article}
\begin{document}
\newcommand{\up}{\vspace*{-0.2cm}}
\newcommand{\upp}{\vspace*{-0.3cm}}
\newcommand{\qed}{\hfill$\rule{.05in}{.1in}$\vspace{.3cm}}
\newcommand{\pf}{\noindent{\bf Proof: }}
\newtheorem{thm}{Theorem}
\newtheorem{lem}{Lemma}
\newtheorem{prop}{Proposition}
\newtheorem{prob}{Problem}
\newtheorem{ex}{Example}
\newtheorem{cor}{Corollary}
\newtheorem{conj}{Conjecture}
\newtheorem{cl}{Claim}
\newtheorem{df}{Definition}
\newtheorem{rem}{Remark}
\newcommand{\beq}{\begin{equation}}
\newcommand{\eeq}{\end{equation}}
\newcommand{\<}[1]{\left\langle{#1}\right\rangle}
\newcommand{\be}{\begin{enumerate}}
\newcommand{\ee}{\end{enumerate}}
\newcommand{\Bul}{\mbox{$\bullet$ } }
\newcommand{\al}{\alpha}
\newcommand{\ep}{\epsilon}
\newcommand{\si}{\sigma}
\newcommand{\om}{\omega}
\newcommand{\la}{\lambda}
\newcommand{\La}{\Lambda}
\newcommand{\Ga}{\Gamma}
\newcommand{\ga}{\gamma}
\newcommand{\im}{\Rightarrow}
\newcommand{\2}{\vspace{.2cm}}
\newcommand{\es}{\emptyset}

\title{\bf Upper Bounds for $\alpha$-Domination Parameters}
\author{A. Gagarin\\
{\footnotesize Department of Mathematics and Statistics}\\
{\footnotesize Acadia University}\\
{\footnotesize Wolfville, Nova Scotia, B4P 2R6}\\
{\footnotesize Canada}\\
\\
A. Poghosyan and V.E. Zverovich\\
{\footnotesize Faculty of Computing, Engineering and Mathematical Sciences}\\
{\footnotesize University of the West of England}\\
{\footnotesize Bristol, BS16 1QY}\\
{\footnotesize UK}}
\date{}
\maketitle

\begin{abstract}
In this paper, we provide a new upper bound for the
$\alpha$-domination number. This result generalises the well-known
Caro-Roditty bound for the domination number of a graph. The same
probabilistic construction is used to generalise another
well-known upper bound for the classical domination in graphs. We
also prove similar upper bounds for the $\alpha$-rate domination number,
which combines the concepts of $\alpha$-domination and $k$-tuple
domination.

\vspace*{.2cm} \noindent {\footnotesize Keywords: {\it Graph;
Domination; $\alpha$-Domination; $\alpha$-Rate Domination}}
\end{abstract}


\section{Introduction}

Domination is one of the fundamental concepts in graph theory with
various applications to ad hoc networks, biological networks,
distributed computing, social networks and web graphs
\cite{alb1,coo1,dai1,hay}. Dominating sets in graphs are natural
models for facility location problems in operational research. An
important role is played by multiple domination, for example
$k$-dominating sets can be used for balancing efficiency and fault
tolerance \cite{dai1}.

We consider undirected simple finite graphs. If $G$ is a graph of
order $n$, then $V(G)=\{v_1,v_2,...,v_n\}$ is the set of vertices
of $G$ and $d_i$ denotes the degree of $v_i$. Let $N(v)$ denote
the neighbourhood of a vertex $v$ in $G$, and $N[v] = N(v)\cup
\{v\}$ be the closed neighbourhood of $v$. A set $X\subseteq V(G)$
is called a {\it dominating set} if every vertex not in $X$ is
adjacent to at least one vertex in $X$. The minimum cardinality of
a dominating set of $G$ is the {\it domination number}
$\gamma(G)$. A set $X$ is called a {\it $k$-dominating set} if
every vertex not in $X$ has at least $k$ neighbors in $X$. The
minimum cardinality of a $k$-dominating set of $G$ is the {\it
$k$-domination number} $\gamma_k(G)$. A set $X$ is called a {\it
$k$-tuple dominating set} of $G$ if for every vertex $v\in V(G)$,
$|N[v]\cap X|\ge k$. The minimum cardinality of a $k$-tuple
dominating set of $G$ is the {\it $k$-tuple domination number}
$\gamma_{\times k}(G)$. The $k$-tuple domination number is only
defined for graphs with $\delta\ge k-1$. A number of upper bounds
for the multiple domination numbers can be found in
\cite{car1,fav1,fav2,GZ,Z2}.

Let $\alpha$ be a real number satisfying $0<\alpha\le1$. A set
$X\subseteq V(G)$ is called an {\it $\alpha$-dominating set} of $G$
if for every vertex $v\in V(G) - X$, $|N(v)\cap X|\ge \alpha d_v$,
i.e. $v$ is adjacent to at least $\lceil\alpha d_v\rceil$ vertices
of $X$. The minimum cardinality of an $\alpha$-dominating set of $G$
is called the {\it $\alpha$-domination number} $\gamma_{\alpha}(G)$.
The $\alpha$-domination was introduced by Dunbar et al. \cite{DHLM}.
It is easy to see that $\gamma(G)\le\gamma_{\alpha}(G)$, and
$\gamma_{\alpha_1}(G)\le\gamma_{\alpha_2}(G)$ for
${\alpha_1}<{\alpha_2}$. Also, $\gamma(G)=\gamma_{\alpha}(G)$ if
$\alpha$ is sufficiently close to $0$.

For an arbitrary graph $G$ with $n$ vertices and $m$ edges, denote
by $\delta=\delta(G)$ and $\Delta=\Delta(G)$ the minimum and
maximum vertex degrees of $G$, respectively. The following results
are proved in \cite{DHLM}:
\begin{equation} \label{Dunbar}
{{\alpha\delta n}\over {\Delta + \alpha\delta}} \le
\gamma_{\alpha}(G) \le {{\Delta n}\over {\Delta +
(1-\alpha)\delta}}
\end{equation}
and
\begin{equation}
{{2\alpha m}\over {(1+\alpha)\Delta}} \le \gamma_{\alpha}(G) \le
{{(2-\alpha)\Delta n - (2-2\alpha)m}\over {(2-\alpha)\Delta}}.
\end{equation}
Interesting results on $\alpha$-domination perfect graphs can be
found in \cite{DRV}. The problem of deciding whether
$\gamma_{\alpha}(G)\le k$ for a positive integer $k$ is known to
be $NP$-complete \cite{DHLM}. Therefore, it is important to have
good upper bounds for the $\alpha$-domination number and efficient
approximation algorithms for finding `small' $\alpha$-dominating
sets.

For $0<\alpha\le 1$, the $\alpha$\emph{-degree} of a graph $G$ is
defined as follows:
$$
{\widehat d}_{\alpha} = {\widehat d}_{\alpha} (G) = {1\over n}
\sum_{i=1}^n \pmatrix{d_i \cr {\lceil \alpha d_i \rceil -1 }}.
$$
In this paper, we use a probabilistic approach to prove that
$$
\ga_{\alpha}(G) \le \left(1-{\widehat{\delta} \over
{(1+\widehat{\delta})}^{1+1/\widehat{\delta}} \, {{\widehat
d}_{\alpha}}^{\;1/\widehat{\delta}}} \right) n,
$$
where $ \widehat{\delta} = \lfloor \delta (1 - \alpha) \rfloor +1
$. This result generalises the well-known upper bound of Caro and
Roditty (\cite{hay}, p.\,48). Using the same probabilistic
construction, we also show that
$$
\ga_{\alpha}(G) \le {\ln(\widehat{\delta}+1) + \ln {\widehat
d}_{\alpha} +1 \over \widehat{\delta}+1} n,
$$
which generalises another well-known upper bound of Alon and
Spencer \cite{alo}, Arnautov \cite{arn}, Lov\'{a}sz \cite{lov} and
Payan \cite{pay}. Finally, we introduce the $\alpha$-rate
domination number, which combines together the concepts of
$\alpha$-domination and $k$-tuple domination, and show that the
$\alpha$-rate domination number satisfies two similar upper
bounds. The random constructions used in this paper also provide
randomized algorithms to find $\alpha$-dominating and
$\alpha$-rate dominating sets satisfying corresponding bounds.

\section{New Upper Bounds for the $\alpha$-Domination Number}

One of the strongest known upper bounds for the domination number is due
to Caro and Roditty:

\begin{thm}[Caro and Roditty \cite{hay}, p.\,48]
\label{Caro} For any graph $G$ with $\delta\ge 1$,
\begin{equation}\label{Car-Rod}
\ga(G) \le \left(1-{\delta \over {(1+\delta)^{1+1/\delta}}}\right)
n.
\end{equation}
\end{thm}

The upper bound (\ref{Car-Rod}) is generalised for the
$\alpha$-domination number in Theorem \ref{main}. Indeed, if $d_i$
are fixed for all $i=1,\ldots,n$, and $\alpha$ is sufficiently
close to 0, then $ \widehat{\delta} = \delta$ (provided $\delta\ge
1$) and ${\widehat d}_{\alpha} =1$.

\begin{thm}\label{main}
 For any graph $G$,
\begin{equation}\label{bound1}
\ga_{\alpha}(G) \le \left(1-{\widehat{\delta} \over
{(1+\widehat{\delta})}^{1+1/\widehat{\delta}} \,{\widehat
d}_{\alpha}^{\;1/\widehat{\delta}}} \right) n,
\end{equation}
where $ \widehat{\delta} = \lfloor \delta (1 - \alpha) \rfloor +1
$.
\end{thm}

\pf Let $A$ be a set formed by an independent choice of vertices
of $G$, where each vertex is selected with the probability
\begin{equation}\label{p_max}
p = 1-{\left(1 \over {{(1 + \widehat{\delta}){\widehat
d}_{\alpha}}}\right)}^{1/\widehat{\delta}}.
\end{equation}
Let us denote $$ B = \{v_i\in V(G) - A\ :\ |N(v_i)\cap A|\le {\lceil
\alpha d_i \rceil -1 }\}.$$ It is obvious that the set $ D = A \cup
B$ is an $\alpha$-dominating set.
 The expectation
of $|D|$ is
\begin{eqnarray*}
E(|D|) &=& E(|A|)+E(|B|)\\
&=& \sum_{i=1}^{n}P(v_i\in A)+ \sum_{i=1}^{n}P(v_i\in B)\\
&=& pn + \sum_{i=1}^{n} \sum_{r=0}^{\lceil \alpha d_i \rceil
-1}\pmatrix{d_i \cr r} p^r {(1-p)}^{d_i -r +1}.
\end{eqnarray*}
It is easy to see that, for $0\le r\le \lceil \alpha d_i
\rceil -1$,
 $$ \pmatrix{d_i \cr r} \le
\pmatrix{d_i \cr {\lceil \alpha d_i \rceil -1 }}\pmatrix{\lceil
\alpha d_i \rceil -1 \cr r}.
$$
Also,
$$
d_i-\lceil \alpha d_i \rceil \ge \lfloor \delta (1-\alpha) \rfloor .
$$
Therefore,
\begin{eqnarray}
\nonumber
E(|D|) &\le& pn + \sum_{i=1}^{n}{\pmatrix{d_i \cr {\lceil \alpha
d_i \rceil -1 }} {(1-p)}^{d_i-\lceil \alpha d_i \rceil +2}}
\sum_{r=0}^{\lceil \alpha d_i \rceil-1}{\pmatrix{\lceil
\alpha d_i \rceil -1 \cr r}p^r (1-p)^{{\lceil \alpha d_i \rceil-1-r}}}\\
\nonumber
&=& pn + \sum_{i=1}^{n}{\pmatrix{d_i \cr {\lceil \alpha
d_i \rceil -1 }} {(1-p)}^{d_i-\lceil \alpha d_i \rceil +2}}\\
\nonumber
&\le& pn + {(1-p)}^{\lfloor \delta (1-\alpha) \rfloor +2}\widehat d_{\alpha}n \\
\label{expr}
&=& pn + {(1-p)}^{\widehat{\delta}+1}\widehat d_{\alpha}n \\
\nonumber &=& \left(1-{\widehat{\delta} \over
{(1+\widehat{\delta})}^{1+1/\widehat{\delta}} \,{\widehat
d}_{\alpha}^{\;1/\widehat{\delta}}} \right) n.
\end{eqnarray}
Note that the value of $p$ in (\ref{p_max}) is chosen to minimize
the expression (\ref{expr}). Since the expectation is an average
value, there exists a particular $\alpha$-dominating set of order
at most $\left(1-{\widehat{\delta} \over
{(1+\widehat{\delta})}^{1+1/\widehat{\delta}} \,{\widehat
d}_{\alpha}^{\;1/\widehat{\delta}}} \right) n$, as required. The
proof of the theorem is complete. \qed

Notice that in some cases Theorem \ref{main} provides a much
better bound than the upper bound in (\ref{Dunbar}). For example,
if $G$ is a 1000-regular graph, then Theorem \ref{main} gives
$\ga_{0.1}(G) < 0.305n$, while (\ref{Dunbar}) yields only
$\ga_{0.1}(G) < 0.527n$.

\begin{cor}\label{Corollary1}
For any graph $G$,
\begin{equation}\label{bound2}
\ga_{\alpha}(G) \le {\ln(\widehat{\delta}+1) + \ln {\widehat
d}_{\alpha} +1 \over \widehat{\delta}+1} n.
\end{equation}
\end{cor}

\pf We put
$$p=\min \left\{ 1, {\ \ln(\widehat{\delta}+1) + \ln{\widehat d}_{\alpha}\over \widehat{\delta}+1} \right\}.$$
Using the inequality $1-p\le e^{-p}$, we can estimate the expression
(\ref{expr}) as follows:
$$
E(|D|) \le pn + e^{-p (\widehat{\delta} +1)}{\widehat d}_{\alpha} n.
$$
If $p=1$, then the result easily follows. If $p={\
\ln(\widehat{\delta}+1) + \ln{\widehat d}_{\alpha}\over
\widehat{\delta}+1}$, then
$$
E(|D|) \le {\ln(\widehat{\delta}+1) + \ln {\widehat d}_{\alpha} +1
\over \widehat{\delta}+1} n,
$$
as required. \qed

Corollary \ref{Corollary1} generalises the following well-known
upper bound independently proved by several authors
\cite{alo,arn,lov,pay}:
\begin{equation}\label{classical}
\ga(G) \le {\ln(\delta+1) + 1 \over \delta+1} n.
\end{equation}

\section{$\alpha$-Rate Domination}

Define a set $X\subseteq V(G)$ to be an {\it $\alpha$-rate
dominating set} of $G$ if for any vertex $v\in V(G)$,
$$
|N[v]\cap X|\ge \alpha d_v.
$$
Let us call the minimum cardinality of an $\alpha$-rate dominating
set of $G$ the {\it $\alpha$-rate domination number}
$\ga_{\times\alpha}(G)$. It is easy to see that
$\gamma_{\alpha}(G)\le \ga_{\times\alpha}(G)$. The concept of
$\alpha$-rate domination is similar to the well-known $k$-tuple
domination (for example, see \cite{KL,Z2}). For $0<\alpha\le 1$,
the \emph{closed $\alpha$-degree} of a graph $G$ is defined as
follows:
$$
{\widetilde d}_{\alpha} = {\widetilde d}_{\alpha} (G) = {1\over n}
\sum_{i=1}^n \pmatrix{d_i+1 \cr {\lceil \alpha d_i \rceil -1 }}.
$$
In fact, the only difference between the $\alpha$-degree and the
closed $\alpha$-degree is that to compute the latter we choose
from $d_i+1$ vertices instead of $d_i$, i.e. from the {\it closed}
neighborhood $N[v_i]$ of $v_i$ instead of $N(v_i)$.

The following theorem provides an analogue of the Caro-Roditty bound
(Theorem \ref{Caro}) for the $\alpha$-rate domination number:

\begin{thm} \label{main3}
For any graph $G$ and $0<\alpha\le1$,
\begin{equation}\label{bound4}
\ga_{\times\alpha}(G) \le \left(1-{\widehat{\delta} \over
{(1+\widehat{\delta})}^{1+1/\widehat{\delta}} \,{\widetilde
d}_{\alpha}^{\;1/\widehat{\delta}}} \right) n,
\end{equation}
where $ \widehat{\delta} = \lfloor \delta (1 - \alpha) \rfloor +1
$.
\end{thm}

\pf Let $A$ be a set formed by an independent choice of vertices of
$G$, where each vertex is selected with probability $p$, $0\le p \le
1$. For $m\ge0$, denote by $B_m$ the set of vertices $v\in V(G)$
dominated by exactly $m$ vertices of $A$ and such that $|N[v]\cap
A|<\alpha d_v$, i.e.
$$
|N[v]\cap A| = m \le \lceil \alpha d_v \rceil - 1.
$$
Note that each vertex $v\in V(G)$ is in at most one of the sets
$B_m$ and $0\le m\le \lceil \alpha d_v \rceil -1$. We form a set
$B$ in the following way: for each vertex $v \in B_m$, select
$\lceil \alpha d_v \rceil - m $ vertices from $N(v)$ that are not
in $A$ and add them to $B$. Consider the set $D=A\cup B$. It is
easy to see that $D$ is an $\alpha$-rate dominating set. The
expectation of $|D|$ is:
\begin{eqnarray*}
E(|D|) &\le& E(|A|) + E(|B|)\\
& \le & \sum_{i=1}^{n}P(v_i\in A)+
\sum_{i=1}^{n}\sum_{m=0}^{\lceil
\alpha d_i \rceil - 1}(\lceil \alpha d_i \rceil - m) P(v_i\in B_m)\\
& = & pn + \sum_{i=1}^{n}\sum_{m=0}^{\lceil \alpha d_i \rceil -
1}(\lceil \alpha d_i \rceil - m) \pmatrix{d_i+1 \cr m} p^m (1-p)^{d_i+1-m}  \\
& \le & pn + \sum_{i=1}^n \sum_{m=0}^{\lceil \alpha d_i \rceil - 1} \pmatrix{d_i+1 \cr \lceil \alpha d_i \rceil-1} \pmatrix{\lceil \alpha d_i \rceil-1 \cr m} p^m (1-p)^{d_i+1-m}\\
& = & pn + \sum_{i=1}^n \pmatrix{d_i+1 \cr \lceil \alpha d_i \rceil-1} (1-p)^{d_i-\lceil \alpha d_i \rceil+2} \sum_{m=0}^{\lceil \alpha d_i \rceil - 1} \pmatrix{\lceil \alpha d_i \rceil-1 \cr m} p^m (1-p)^{\lceil \alpha d_i \rceil-1-m}\\
& = & pn + \sum_{i=1}^n \pmatrix{d_i+1 \cr \lceil \alpha d_i \rceil-1} (1-p)^{d_i-\lceil \alpha d_i \rceil+2}\\
& \le & pn +  (1-p)^{\lfloor \delta (1-\alpha) \rfloor +2}
\sum_{i=1}^n
\pmatrix{d_i+1 \cr \lceil \alpha d_i \rceil-1}\\
& = & pn + (1-p)^{\widehat{\delta}+1} {\widetilde d}_{\alpha} n,
\end{eqnarray*}
since
$$(\lceil \alpha d_i \rceil - m) \pmatrix{d_i+1 \cr m} \le \pmatrix{d_i+1 \cr \lceil \alpha d_i \rceil-1} \pmatrix{\lceil \alpha d_i \rceil-1 \cr m}.$$
Thus,
\begin{equation}\label{again}
E(|D|) \le pn + (1-p)^{\widehat{\delta}+1} {\widetilde d}_{\alpha}
n.
\end{equation}
Minimizing the expression (\ref{again}) with respect to $p$, we
obtain
$$
E(|D|) \le \left(1-{\widehat{\delta} \over
{(1+\widehat{\delta})}^{1+1/\widehat{\delta}} \,{\widetilde
d}_{\alpha}^{\;1/\widehat{\delta}}} \right) n,
$$
as required. The proof of Theorem \ref{main3} is complete. \qed

\begin{cor}\label{Corollary2}
For any graph $G$,
\begin{equation}\label{bound5}
\ga_{\times\alpha}(G) \le {\ln(\widehat{\delta}+1) + \ln {\widetilde
d}_{\alpha} +1 \over \widehat{\delta}+1} n.
\end{equation}
\end{cor}

\pf Using an approach similar to that in the proof of Corollary
\ref{Corollary1}, the result follows if we put
$$p=\min \left\{ 1, {\ \ln(\widehat{\delta}+1) + \ln{\widetilde d}_{\alpha}\over \widehat{\delta}+1} \right\}$$
and use the inequality $1-p\le e^{-p}$ to estimate the expression
(\ref{again}) as follows:
$$ E(|D|) \le pn + e^{-p (\widehat{\delta} +1)}{\widetilde d}_{\alpha} n.
$$
\qed

Note that, similar to Corollary \ref{Corollary1}, the bound of
Corollary \ref{Corollary2} also generalises the classical upper
bound (\ref{classical}). However, the probabilistic construction
used to obtain the bounds (\ref{bound4}) and (\ref{bound5}) is
different from that to obtain the bounds (\ref{bound1}) and
(\ref{bound2}).

\section{Final Remarks and Open Problems}

Notice that the concept of the $\alpha$-rate domination number
$\ga_{\times\alpha}(G)$ is `opposite' to the $\alpha$-independent
$\alpha$-domination number $i_{\alpha}(G)$ as defined in
\cite{DRV}. It would be interesting to use a probabilistic method
construction to obtain an upper bound for $i_{\alpha}(G)$.

Also, the random constructions used to obtain the upper bounds
(\ref{bound1}), (\ref{bound2}), (\ref{bound4}) and (\ref{bound5})
provide randomized algorithms to find corresponding dominating
sets in a given graph $G$. It would be interesting to derandomize
these algorithms or to obtain independent deterministic algorithms
to find corresponding dominating sets satisfying the upper bounds
(\ref{bound1}), (\ref{bound2}), (\ref{bound4}) and (\ref{bound5}).
Algorithms approximating the $\alpha$- and $\alpha$-rate
domination numbers up to a certain degree of precision would be
interesting too. For the $k$-tuple domination number, an
interesting approximation algorithm was found by Klasing and
Laforest \cite{KL}.

Using probabilistic methods, Alon \cite{alo2} proved that the
bound (\ref{classical}) is asymptotically best possible. More
precisely, it was proved that when $n$ is large there exists a
graph $G$ such that
$$
\ga(G) \ge {\ln(\delta+1)+1 \over \delta+1} n (1+o(1)).
$$
We wonder if a similar result can be proved for the bounds
(\ref{bound2}) and (\ref{bound5}), and conjecture that when $n$ is
large enough there exist graphs $G$ and $H$ such that
$$\ga_{\alpha}(G) \ge {\ln(\widehat{\delta}+1) + \ln {\widehat
d}_{\alpha} +1 \over \widehat{\delta}+1}n(1+o(1))
$$
and
$$
\ga_{\times\alpha}(H) \ge {\ln(\widehat{\delta}+1) + \ln
{\widetilde d}_{\alpha} +1 \over \widehat{\delta}+1} n(1+o(1)).
$$

\end{document}